\documentclass[twoside]{amsart}

\usepackage{amsmath,amssymb,amsbsy,amstext,amsxtra,latexsym}
\usepackage{graphicx}
\usepackage{enumerate}
\usepackage[all]{xy}

\newtheorem{theorem}{Theorem}

\newtheorem{proposition}{Proposition}

\theoremstyle{remark}

\newcommand{\uc}[1]{\ensuremath \overset{#1}{\circ}}

\def\mP{{\mathbb P}}
\def\mQ{{\mathbb Q}}

\def\mZ{{\mathbb Z}}

\begin{document}

\title[A complex surface of general type with $p_g=0$, $K^2=3$ and
       $H_1=\mZ/2\mZ$]{A complex surface of general type with
       $p_g=0$, $K^2=3$ and $H_1=\mZ/2\mZ$}

\author{Heesang Park, Jongil Park and Dongsoo Shin}

\address{Department of Mathematical Sciences, Seoul
         National University,
         San 56-1, Sillim-dong, Gwanak-gu, Seoul 151-747, Korea}

\email{hspark@math.snu.ac.kr}

\address{Department of Mathematical Sciences, Seoul
         National University,
         San 56-1, Sillim-dong, Gwanak-gu, Seoul 151-747, Korea}

\email{jipark@math.snu.ac.kr}

\address{Department of Mathematical Sciences, Seoul
         National University,
         San 56-1, Sillim-dong, Gwanak-gu, Seoul 151-747, Korea}

\email{dsshin@math.snu.ac.kr}

\date{March 13, 2008}

\subjclass[2000]{Primary 14J29; Secondary 14J10, 14J17, 53D05}

\keywords{$\mQ$-Gorenstein smoothing, rational blow-down, surface of general type}

\maketitle

 This paper is an addendum to~\cite{PPS}, in which the authors constructed
 a simply connected minimal complex surface of general type with $p_g=0$
 and $K^2=3$.
 Motivated by Y. Lee and the second author's recent construction on
 a surface of general type with $p_g=0$, $K^2=2$ and $H_1=\mZ/2\mZ$~\cite{LP08},
 we extend the result to the $K^2=3$ case in this paper.
 That is, we construct a new non-simply connected minimal
 surface of general type with $p_g=0$, $K^2=3$ and $H_1=\mZ/2\mZ$
 using a rational blow-down surgery and a $\mQ$-Gorenstein smoothing theory.

 The key ingredient of this paper is to find a right rational surface $Z$
 which makes it possible to get such a complex surface.
 Once we have a right candidate $Z$ for $K^2=3$, the remaining argument
 is similar to that of $K^2=3$ case appeared in~\cite{PPS}.
 That is, by applying a rational blow-down surgery and a $\mQ$-Gorenstein
 smoothing theory developed in~\cite{LP07} to $Z$,
 we obtain a minimal complex surface of general type with $p_g=0$ and $K^2 =3$.
 Then we show that the surface has $H_1 = \mZ/2\mZ$
 by using a similar method in~\cite{LP08}.
 Since almost all the proofs are parallel to the case of the main
 construction in~\cite[\S3]{PPS},
 we only explain how to construct such a minimal complex surface.
 The main result of this paper is the following

\begin{theorem}
\label{thm-main}
 There exists a minimal complex surface of general type with $p_g=0$,
 $K^2=3$ and $H_1=\mZ/2\mZ$.
\end{theorem}

\subsubsection*{Remark}
 Although some examples of non-simply connected
 complex surfaces of general type with $p_g =0$ and $K^2=3$ have been
 constructed (\cite{BHPV}, VII), until now it is not known whether there is
 a complex surface of general type with $p_g=0$, $K^2 =3$ and $H_1=\mZ/2\mZ$.
 Theorem~\ref{thm-main} above provides the first example of such a complex
 surface.

\section{Main construction}

\begin{figure}[hbtb]
 \centering
 \includegraphics[scale=0.7]{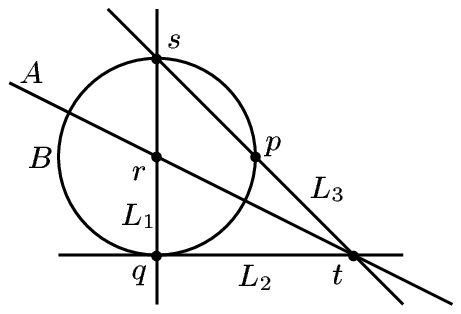}
 \caption{A pencil of cubics}
 \label{figure:pencil}
\end{figure}

\begin{figure}[hbtb]
\centering
\includegraphics[scale=0.7]{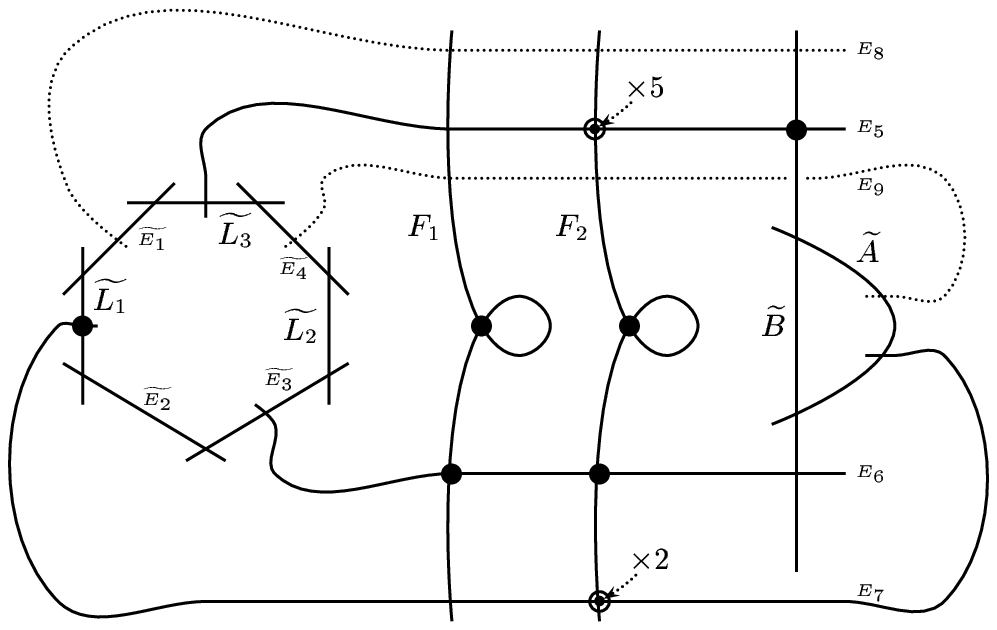}
\caption{An elliptic fibration $Y$}
\label{figure:Y}
\end{figure}

 We start with a special elliptic fibration $Y:=\mP^2\sharp 9\overline{\mP}^2$
 which is used in the main construction of this paper.
 Let $L_1$, $L_2$, $L_3$ and $A$ be lines in $\mP^2$ and let $B$ be a smooth
 conic in $\mP^2$ intersecting as in Figure~\ref{figure:pencil}.
 We consider a pencil of cubics
 $\{ \lambda(L_1+L_2+L_3) + \mu(A+B) \mid [\lambda:\mu] \in \mP^1 \}$
 in $\mP^2$ generated by two cubic curves $L_1+L_2+L_3$ and $A+B$,
 which has $5$ base points, say, $p$, $q$, $r$, $s$ and $t$.
 In order to obtain an elliptic fibration over $\mP^1$ from the pencil,
 we blow up three times at $q$ and twice at $s$ and $t$, respectively,
 including infinitely near base-points at each point.
 We perform two further blowing-ups at the base points $p$ and $r$.
 By blowing-up nine times, we resolve all base points
 (including infinitely near base-points) of the pencil
 and we then get an elliptic fibration $Y=\mP^2\sharp 9\overline{\mP}^2$
 over $\mP^1$ (Figure~\ref{figure:Y}).

 We denote by $E_i$ (or $\widetilde{E_i}$), $i=1,\dotsc,9$, the exceptional divisors (or their proper transforms in $Y$, respectively) induced by the nine blowing-ups.
 Note that there are five sections of the elliptic
 fibration $Y$ corresponding to the five base points $p$, $q$, $r$, $s$,
 and $t$, which may be denoted by $E_5, E_6, \dotsc, E_9$, respectively.
 Furthermore, the elliptic fibration $Y$ has an $I_7$-singular fiber
 consisting of the proper transforms $\widetilde{L_i}$ of $L_i$ ($i=1,2,3$),
 $\widetilde{E_1}$, $\widetilde{E_2}$, $\widetilde{E_3}$ and $\widetilde{E_4}$.
 Also $Y$ has an $I_2$-singular fiber consisting of the proper transforms
 $\widetilde{A}$ and $\widetilde{B}$ of $A$ and $B$, respectively.
 According to the list of Persson~\cite{Pers},
 we may assume that $Y$ has three more nodal singular fibers by choosing
 generally $L_i$'s, $A$ and $B$. Among the three nodal singular fibers,
 we use only two nodal singular fibers, say $F_1$ and $F_2$, for the main
 construction (Figure~\ref{figure:Y}).

\begin{figure}[hbtb]
\centering
\includegraphics[scale=0.7]{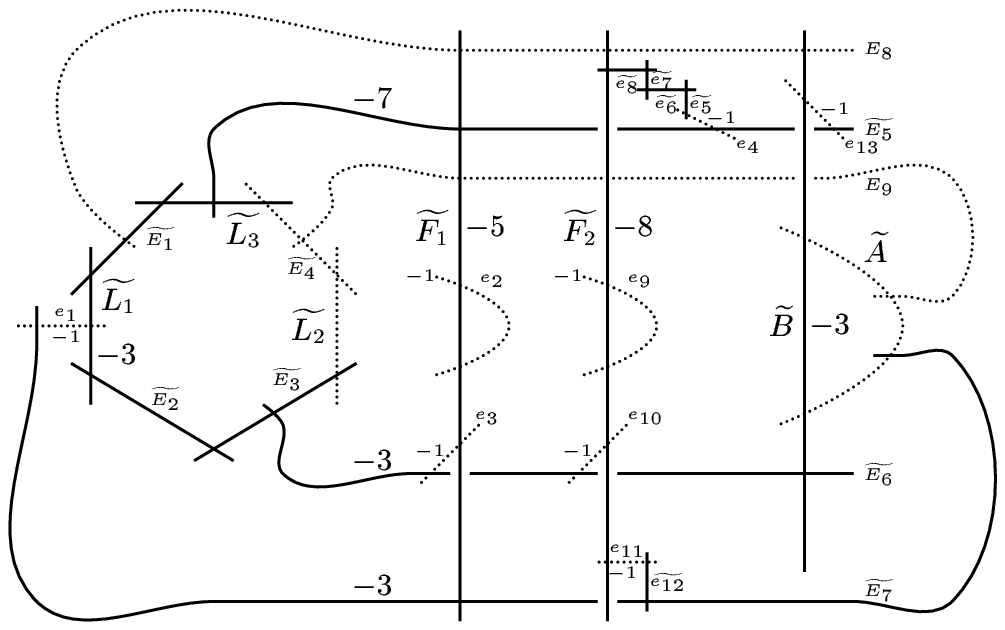}
\caption{A rational surface $Z=Y\sharp 13{\overline \mP}^2$}
\label{figure:Z}
\end{figure}

 Next, by blowing-up several times on $Y$, we construct a rational surface
 $Z$ which contains special configurations of linear chains of $\mP^1$'s.
 At first we blow up five times at the marked point $\bigodot$ on
 $F_2 \cap E_5$. We also blow up two times at the marked point
 $\bigodot$ on $F_2 \cap E_7$.
 Finally we blow up at the six marked points $\bullet$ on each fiber.
 We then get a rational surface $Z=Y \sharp 13\overline{\mP}^2$.
 We denote by $e_i$ (or $\widetilde{e_i}$), $i=1,\dotsc,13$, the exceptional divisors (or their
 proper transforms in $Z$, respectively) induced by the $13$ blow-ups and
 we also denote by $\widetilde{F_i}$ ($i=1,2$) the proper transforms of $F_i$.
 Then there exist two disjoint linear chains of ${\mP}^1$'s in $Z$:
 $C_{110,67} = \uc{-2}-\uc{-3}-\uc{-5}-\uc{-7}-\uc{-2}-\uc{-2}-\uc{-3}-
 \uc{-2}-\uc{-2}-\uc{-3}-\uc{-3}$
 (which consists of $\widetilde{e_{12}}$, $\widetilde{E_7}$, $\widetilde{F_1}$, $\widetilde{E_5}$, $\widetilde{L_3}$, $\widetilde{E_1}$, $\widetilde{L_1}$, $\widetilde{E_2}$, $\widetilde{E_3}$, $\widetilde{E_6}$, $\widetilde{B}$) and $C_{6,1} = \uc{-8}-\uc{-2}-\uc{-2}-\uc{-2}-\uc{-2}$ (which consists of
 $\widetilde{F_2}$, $\widetilde{e_8}$, $\widetilde{e_7}$, $\widetilde{e_6}$, $\widetilde{e_5}$) (Figure~\ref{figure:Z}).

 Next, by applying $\mQ$-Gorenstein smoothing theory as in~\cite{PPS},
 we construct the minimal complex surface appeared in the main theorem.
 That is, we first contract two disjoint chains $C_{110,67}$ and $C_{6,1}$
 of $\mP^1$'s from $Z$ so that it produces a normal projective surface $X$
 with two permissible singular points.
 And then, by using a similar technique in~\cite{PPS},
 we can conclude that $X$ has a $\mQ$-Gorenstein smoothing
 and a general fiber $X_t$ of the $\mQ$-Gorenstein smoothing of $X$
 is a minimal complex surface of general type with $p_g=0$ and $K^2=3$.
 Let us denote a general fiber of the $\mQ$-Gorenstein smoothing
 of $X$ by $X_t$.
 Finally it remains to show that $H_1(X_t;\mZ)=\mZ/2\mZ$.

\medskip

\subsection*{Proof of $H_1(X_t;\mZ)=\mZ/2\mZ$}

 Let $Z_{110,6}$ be a rational blow-down $4$-manifold obtained from $Z$
 by replacing two disjoint configurations $C_{110,67}$ and $C_{6,1}$
 with the corresponding rational balls $B_{110,67}$ and $B_{6,1}$,
 respectively. Then, since a general fiber $X_t$ of a $\mQ$-Gorenstein
 smoothing of $X$ is diffeomorphic to the rational blow-down $4$-manifold
 $Z_{110,6}$, we have $H_1(X_t;\mZ)=H_1(Z_{110,6};\mZ)$.
 Hence it suffices to show that $H_1(Z_{110,6};\mZ) = \mZ/2\mZ$.

\begin{proposition}
\label{prop-0orZ2}
 $H_1(Z_{110,6};\mZ) = \mZ/2\mZ$.
\end{proposition}

\begin{proof}
 First note that the rational surface $Z=Y \sharp 13\overline{\mP}^2$
 can be decomposed into $Z = Z_0 \cup \{ C_{110,67} \cup C_{6,1}\}$
 and the rational blow-down $4$-manifold $Z_{110,6}$ can be
 decomposed into $Z_{110,6} = Z_0 \cup \{B_{110,67} \cup B_{6,1}\}$.

 Let $W = Z_0 \cup B_{110,67}$ and consider the following exact homology
 sequence for a pair $(W, \partial W)$:
\begin{equation*}
 \cdots \rightarrow H_2(W, \partial W;\mZ) \xrightarrow{\partial_{\ast}}
 H_1(\partial W;\mZ) \xrightarrow{i_{\ast}} H_1(W;\mZ) \rightarrow 0.
\end{equation*}
 Note that $\partial W = L(36, -5)$ and a generator of
 $H_1(\partial W;\mZ) = \mZ/36\mZ$ can be represented by a normal circle,
 say $\alpha$, of a disk bundle $C_{6,1}$ over $(-8)$-curve $\widetilde{F_2}$.
 Then we have
\begin{equation*}
 \partial_{\ast}([e_9|_W]) = 2\alpha \in  H_1(\partial W;\mZ) = \mZ/36\mZ.
\end{equation*}
 Furthermore, by using a similar technique in Section 2 of~\cite{LP08},
 we can conclude that the generator $\alpha \in H_1(\partial W;\mZ)$
 is not in the image of $\partial_{\ast}$. Hence it follows from the
 exact sequence above that we have $H_1(W;\mZ) = \mZ/2\mZ$ generated
 by an element $i_{\ast}(\alpha)$.

 Next, we consider the Mayer-Vietoris sequence for a triple
 $(Z_{110,6};W,B_{6,1})$:
\begin{equation*}%\scriptstyle
 H_2(Z_{110,6};\mZ) \xrightarrow{\partial_{\ast}} H_1(L(36,-5);\mZ)
 \xrightarrow{i_{\ast} \oplus j_{\ast}} H_1(W;\mZ) \oplus
 H_1(B_{6,1};\mZ) \rightarrow H_1(Z_{110,6};\mZ) \rightarrow 0.
\end{equation*}
 Since $i_{\ast} \oplus j_{\ast} : H_1(L(36,-5);\mZ)
 \to H_1(W;\mZ) \oplus H_1(B_{6,1};\mZ)$ sends a generator,
 $\alpha$, to (a generator, a generator),
 we finally have $H_1(Z_{110,6};\mZ) = \mZ/2\mZ$.
\end{proof}

\bigskip

\subsubsection*{Acknowledgements}
 Jongil Park was supported by the Korea Research Foundation Grant funded
 by the Korean Government (KRF-2007-314-C00024) and he also holds a
 joint appointment in the Research Institute of Mathematics, SNU.
 Dongsoo Shin was supported by Korea Research Foundation Grant funded
 by the Korean Government (KRF-2005-070-C00005).

\bigskip
\bigskip

\providecommand{\bysame}{\leavevmode\hbox
 to3em{\hrulefill}\thinspace}

\end{document}